\documentclass[11pt]{article}
\usepackage{amsfonts,amsbsy,amsmath,amssymb}
\usepackage{graphics}
\usepackage[lmargin=3.0cm, rmargin=3.0cm, tmargin=3.0cm, bmargin=3.0cm]{geometry}

\numberwithin{equation}{section}

\newcommand{\iso}{\cong} 

\newcommand{\scaps}[1]{{\scshape #1}}
\newcommand{\bfm}{\textbf}
\newcommand{\mcal}{\mathcal}

\newcommand{\raf}[1]{(\ref{#1})}

\renewcommand{\P}{\mcal{P}}

\newcommand{\bnd}[2]{{\rm bnd}_{#2}#1}

\newcommand{\sm}{\setminus}

\def\QED{$\blacksquare$}
\def\inQED{$\square$}

\newenvironment{statement}
{
\refstepcounter{equation} 
\ \\
\noindent
\begin{it}
\noindent
\emph{\bfm{~\theequation.}}
}
{\end{it}\ \\}

\newcommand{\sref}[1]{\bfm{{\ref{#1}}}}

\makeatletter
\def\section{\@ifstar\unnumberedsection\numberedsection}
\def\numberedsection{\@ifnextchar[
  \numberedsectionwithtwoarguments\numberedsectionwithoneargument}
\def\unnumberedsection{\@ifnextchar[
  \unnumberedsectionwithtwoarguments\unnumberedsectionwithoneargument}
\def\numberedsectionwithoneargument#1{\numberedsectionwithtwoarguments[#1]{#1}}
\def\unnumberedsectionwithoneargument#1{\unnumberedsectionwithtwoarguments[#1]{#1}}
\def\numberedsectionwithtwoarguments[#1]#2{%
  \ifhmode\par\fi
  \removelastskip
  \vskip 3ex\goodbreak
  \refstepcounter{section}%
  \noindent
  \leavevmode
  \begingroup
  \bfseries
  \S \thesection\ 
  #2.\quad
  \endgroup
  \addcontentsline{toc}{section}{%
    \protect\numberline{\thesection}%
    #1}%
  }
\def\unnumberedsectionwithtwoarguments[#1]#2{%
  \ifhmode\par\fi
  \removelastskip
  \vskip 3ex\goodbreak
  \noindent
  \leavevmode
  \begingroup
  \bfseries
  #2.\quad
  \endgroup
  \addcontentsline{toc}{section}{%
    #1}%
  }

\makeatletter
\def\subsection{\@ifstar\unnumberedsubsection\numberedsubsection}
\def\numberedsubsection{\@ifnextchar[
  \numberedsubsectionwithtwoarguments\numberedsubsectionwithoneargument}
\def\unnumberedsubsection{\@ifnextchar[
  \unnumberedsubsectionwithtwoarguments\unnumberedsubsectionwithoneargument}
\def\numberedsubsectionwithoneargument#1{\numberedsubsectionwithtwoarguments[#1]{#1}}
\def\unnumberedsubsectionwithoneargument#1{\unnumberedsubsectionwithtwoarguments[#1]{#1}}
\def\numberedsubsectionwithtwoarguments[#1]#2{%
  \ifhmode\par\fi
  \removelastskip
  \vskip 3ex\goodbreak
  \refstepcounter{subsection}%
  \noindent
  \leavevmode
  \begingroup
  \bfseries
  \S \thesubsection\ 
  #2.\quad
  \endgroup
  \addcontentsline{toc}{subsection}{%
    \protect\numberline{\thesubsection}%
    #1}%
  }
\def\unnumberedsubsectionwithtwoarguments[#1]#2{%
  \ifhmode\par\fi
  \removelastskip
  \vskip 3ex\goodbreak
  \noindent
  \leavevmode
  \begingroup
  \bfseries
  #2.\quad
  \endgroup
  \addcontentsline{toc}{subsection}{%
    #1}%
  }
\makeatother

\makeatletter
\def\subsubsection{\@ifstar\unnumberedsubsubsection\numberedsubsubsection}
\def\numberedsubsubsection{\@ifnextchar[
  \numberedsubsubsectionwithtwoarguments\numberedsubsubsectionwithoneargument}
\def\unnumberedsubsubsection{\@ifnextchar[
  \unnumberedsubsubsectionwithtwoarguments\unnumberedsubsubsectionwithoneargument}
\def\numberedsubsubsectionwithoneargument#1{\numberedsubsubsectionwithtwoarguments[#1]{#1}}
\def\unnumberedsubsubsectionwithoneargument#1{\unnumberedsubsubsectionwithtwoarguments[#1]{#1}}
\def\numberedsubsubsectionwithtwoarguments[#1]#2{%
  \ifhmode\par\fi
  \removelastskip
  \vskip 3ex\goodbreak
  \refstepcounter{subsubsection}%
  \noindent
  \leavevmode
  \begingroup
  \bfseries
  \S \thesubsubsection\ 
  #2.\quad
  \endgroup
  \addcontentsline{toc}{subsubsection}{%
    \protect\numberline{\thesubsubsection}%
    #1}%
  }
\def\unnumberedsubsubsectionwithtwoarguments[#1]#2{%
  \ifhmode\par\fi
  \removelastskip
  \vskip 3ex\goodbreak
  \noindent
  \leavevmode
  \begingroup
  \bfseries
  #2.\quad
  \endgroup
  \addcontentsline{toc}{subsubsection}{%
    #1}%
  }
\makeatother

\begin{document}
\begin{center}
{\Large \scaps{The Kelmans-Seymour conjecture for apex graphs}}\\ \ \\

{\small
Elad Aigner-Horev\footnote{horevel@cs.bgu.ac.il.} and
Roi Krakovski\footnote{roikr@cs.bgu.ac.il.} \\ \ \\  

Department of Computer Science\\
Ben-Gurion University of the Negev,\\ 
Beer Sheva, 84105, Israel}
\end{center}

{\small
\noindent
\bfm{Abstract.} We provide a short proof that a $5$-connected nonplanar apex graph contains a subdivided $K_{_5}$ or a $K^-_{_4}$ (= $K_{_4}$ with a single edge removed) as a subgraph. Together with a recent result of Ma and Yu that {\sl every nonplanar $5$-connected graph containing $K^-_{_4}$ as a subgraph has a subdivided $K_{_5}$}; this settles the Kelmans-Seymour conjecture for apex graphs.\\ 

\noindent
\scaps{Keywords.} Subdivided $K_{_5}$, Apex graphs. \\

\noindent 
\scaps{Preamble.} Whenever possible notation and terminology are that of~\cite{diestel}.
Throughout, a graph is always simple, undirected, and finite. $G$ always denotes a graph. A subdivided $G$ is denoted $TG$. $K^-_{_4}$ denotes $K_{_4}$ with a single edge removed.
We write $\delta(G)$ and $d_G(v)$ to denote the minimum degree of $G$ and the degree of a vertex $v \in V(G)$, respectively. The $k$-\emph{wheel} graph consists of a $k$-circuit $C$ and an additional vertex, called the \emph{hub}, adjacent to every vertex of $C$ through edges called the \emph{spokes}. $C$ is called the \emph{rim} of the wheel.} 

\section{Introduction} A refinement of Kuratowski's theorem postulated by the Kelmans-Seymour conjecture ($1975$) is that: {\sl the $5$-connected nonplanar graphs contain a $TK_{_5}$}. As this conjecture is open for many years now, it does not stand to reason that certain special cases of this conjecture be considered. If to pick a special case, then we contend that the apex graphs are a natural choice; where a graph is \emph{apex} if it has a vertex, referred to as an \emph{apex vertex}, removal of which results in a planar graph.  
In this paper, we prove in a short manner that:

\begin{statement}\label{main}
A $5$-connected nonplanar apex graph contains a $TK_{_5}$ or a $K^-_{_4}$ as a subgraph. 
\end{statement}

Recently, Ma and Yu~\cite{MaYu1,MaYu2} proved that:

\begin{statement}\label{MaYuk4}
\bfm{(Ma-Yu~\cite{MaYu1,MaYu2})}\\
A $5$-connected nonplanar graph containing $K^-_{_4}$ as a subgraph, contains a $TK_{_5}$. 
\end{statement}

\noindent
By \sref{main} and \sref{MaYuk4}, it follows that the Kelmans-Seymour conjecture holds for apex graphs. 

\begin{statement}
A $5$-connected nonplanar apex graph contains a $TK_{_5}$. 
\end{statement}

\noindent
\scaps{Our proof of} \sref{main}. By Euler's formula, a $2$-connected planar graph with minimum degree $5$ contains $K^-_{_4}$ as a subgraph~\cite[Lemma $2$]{Mohar}. 
Consequently, a $5$-connected nonplanar apex graph $G$ satisfying $K^-_{_4} \not\subseteq G$ has $\delta(G-v) = 4$, where $v$ is an apex vertex of $G$. Thus, a $5$-connected nonplanar apex graph contains $K^-_{_4}$ as a subgraph or has an apex vertex that is part of a $5$-(vertex)-disconnector of $G$. Thus, to prove \sref{main}, suffices that we prove the following. 

\begin{statement}\label{main2}
A $5$-connected nonplanar apex graph $G$ with an apex vertex contained in a $5$-(vertex)-disconnector of $G$ satisfies $TK_{_5} \subseteq G$ or $K^-_{_4} \subseteq G$.  
\end{statement} 

Adjourning technical details until later sections, we outline here the sole manner in which we construct a $TK_{_5}$ in our proof of \sref{main2} assuming $K^-_{_4} \not\subseteq G$ and $v$ is an apex vertex of $G$ satisfying the premise of \sref{main2}.
\begin{enumerate}
\item [(S.1)] We fix an embedding of $G-v$ and identify it with its embedding. We then pick a ``suitable'' $5$-(vertex)-disconnecter $D$ containing $v$ such that 
      $G= G_{_1} \cup G_{_2}$ and $G[D] = G_{_1} \cap G_{_2}$.
\item [(S.2)] In one of the sides of this disconnector, say $G_{_1}$, we find a $4$-valent vertex $u$
      such that together with $u$ the vertices cofacial with $u$ in $G-v$ 
      induce a subdivided $d(u)$-wheel $S \subseteq G_{_1}$ whose spokes are
      preserved and coincide with the edges incident with $u$.
\item [(S.3)] In $G_{_1} -v$, we construct $3$ pairwise vertex-disjoint paths (i.e., a $3$-linkage)
      linking $D-v$ and the rim of $S$ (not meeting $u$) so that these paths meet the rim of 
      $S$ only at $N_G(u)$. 
\item [(S.4)] We choose an arbitrary vertex in $G_{_2}-D$ and connect it to $D$ through 
      a $5$-fan contained in $G_{_2}$. 
\item [(S.5)] $uv \in E(G)$ as $u$ is $4$-valent.
\item [(S.6)] $TK_{_5} \subseteq$ the union of $S$, the $3$-linkage, the $5$-fan, and $uv$. 
\end{enumerate}
Essentially, the remainder of this paper consists of our preparation for this single construction. The accurate form of this construction can be found in \S\ref{proof}. We use the discharging method for finding the wheel $S$ in (S.2). 

\section{Preliminaries}\label{pre}

\noindent
\scaps{Subgraphs.} Let $H$ be a subgraph of $G$, denoted $H \subseteq G$.
The boundary of $H$, denoted by $\bnd{H}{G}$ (or simply $\bnd{H}{}$), is the set of vertices of $H$ incident with $E(G)\sm E(H)$. By $int_G H$ (or simply $int H$) we denote the subgraph induced by $V(H) \sm bnd H$. If $v \in V(G)$, then $N_H(v)$ denotes $N_G(v) \cap V(H)$.\\

\noindent
\scaps{Paths and circuits.} 
For $X,Y \subseteq V(G)$, an $(X,Y)$-path is a simple path with one end in $X$ and the other in $Y$ internally-disjoint of $X \cup Y$. If $X =\{x\}$, we write $(x,Y)$-path.
If $|X|=|Y|=k \geq 1$, then a set of $k$ pairwise vertex-disjoint $(X,Y)$-paths is called an $(X,Y)$-$k$-\emph{linkage}. Throughout this paper, a linkage is always of size $4$. 

If $x \in V(G)$ and $Y \subseteq V(G) \sm \{x\}$, then by $(x,Y)$-$k$-\emph{fan} we mean a set of $k \geq 1$ $(x,Y)$-paths with only $x$ as a common vertex.

The \emph{interior} of an $xy$-path $P$ is the set $V(P)\sm\{x,y\}$ and is denoted $int P$. For $u,v \in V(P)$, we write $[uPv]$ to denote the $uv$-subpath of $P$. We write $(uPv)$ to denote $int [uPv]$, and in a similar manner the semi-open segments $[uPv)$ and $(uPv]$.

If $C$ is a circuit of a plane graph $G$ and $A=\{a_{_1},a_{_2},a_{_3},a_{_4}\} \subseteq V(C)$ appear in this clockwise order along $C$, then $[a_{_i}Ca_{_{i+1}}]$, $1 \leq i \leq 4$, denotes the segment of $C$ whose ends are $a_{_i}$ and $a_{_{i+1}}$ and such that its interior, denoted $(a_{_i}Ca_{_{i+1}})$, does not meet $A$ (clearly, $a_{_5} = a_{_1}$). Semi-open segments $[a_{_i}C_{_u}a_{_{i+1}})$ are defined accordingly.
Two members of $A$ are called \emph{consecutive} if these are consecutive in the clockwise ordering of $A$ along $C$.\\

\noindent
\scaps{Bridges.} Let $H \subseteq G$. By $H$-\emph{bridge} we mean either an edge $uv \notin E(H)$ and $u,v \in V(H)$ or a connected component of $G-H$. In the latter case, the $H$-bridge is called \emph{nontrivial}. The vertices of $H$ adjacent to an $H$-bridge $B$ are called the \emph{attachment vertices} of $B$. A $uv$-path internally-disjoint of $H$ with $u,v \in V(H)$, is called an $H$-\emph{ear}.\\ 

\noindent
\scaps{Hammocks.} A $k$-\emph{hammock} of $G$ is a connected subgraph $H$ satisfying $|bnd H| = k \geq 1$. A hammock $H$ coinciding with its boundary is called {\sl trivial}, \emph{degenerate} if $|V(H)| = |bnd H| +1$, and \emph{fat} if $|V(H)|\geq |bnd H|+2$. We call a $4$-hammock {\sl minimal} if all its proper $4$-hammocks, if any, are trivial or degenerate.

\begin{statement}\label{2-con-ham}
A minimal fat $4$-hammock $H$ of a $4$-connected graph $G$, $K^-_{_4} \not\subseteq G$, satisfies $\kappa(H) \geq 2$.
\end{statement}

\noindent
{\sl Proof.} Assume, to the contrary, that $H = H_{_1} \cup H_{_2}$ such that $\{x\} = V(H_{_1}) \cap V(H_{_2})$ and $V(H_{_i}) \sm \{x\} \not= \emptyset$, for $i=1,2$. Clearly, $bnd H_{_i} = \{x\} \cup X_{_i}$, where $X_{_i} \subset bnd H$, for $i=1,2$. 

Consequently, if $H_{_i}$ is a $4$-hammock of $G$, then $H_{_{3-i}}$ consists of a single edge; implying that $H_{_i}$ is degenerate, by fatness of $H$. As $d_G(x) \geq 4$, $K^-_{_4} \subseteq G$. 

Next, if each of $H_{_i}$, $i=1,2$, is a $k$-hammock of $G$ with $k \leq 3$, then
both are trivial, by $4$-connectivity of $G$. This in turn implies that $H$ is degenerate satisfying $\{x\} = V(H) \sm bnd H$; contradiction to the fatness of $H$.\QED \\

\noindent
\scaps{Subdivided wheels.} For $u \in V(G)$, we write $S_{_u}$ to denote a subdivided $d(u)$-wheel with hub $u$, the spokes preserved and coinciding with $\{uv:v \in N(u)\}$. Its rim, denoted $C_{_u}$, is an induced circuit of $G$ separating $u$ from the rest of $G$. 

If $G$ is a $4$-connected plane graph, then such an $S_{_u}$ exists for every $u \in V(G) \sm V(X_{_G})$, where $X_{_G}$ is the infinite face of $G$. Indeed, 
the set of vertices cofacial with $u$ form $C_{_u}$. Consequently, if $G$ is a plane graph and $u \in V(G) \sm V(X_{_G})$ we refer to $S_{_u}$ as the \emph{facial wheel of $u$}. Such a subdivided wheel is called \emph{short} if:\\
(SH.1) $d(u) =4$ and $u$ is the common vertex of two edge disjoint triangles, say $T$ and $T'$; and\\
(SH.2) the two segments of $C_{_u} - (E(C_{_u}) \cap E(T)) - (E(C_{_u}) \cap E(T))$, say $Q$ and $Q'$, satisfy:\\
\indent (SH.2.a) $2 \leq |V(Q)|,|V(Q')| \leq 4$; and\\
\indent (SH.2.b) if one segment is of order $4$, then the other is of order $\leq 3$.\\

A short wheel is called \emph{imbalanced} if one of its segments is of order $4$.
An imbalanced wheel $S_{_u} \subseteq H$, where $H$ is a $4$-hammock of a $4$-connected graph, is called \emph{proper with respect to $H$} if the interior of its segment of order $4$ does not meet $bnd H$. If $H$ is understood, then we write \emph{proper}.\\ 

\noindent
\scaps{Faces of plane graphs.} Let $G$ be a 2-connected plane graph. By $F(G)$ we denote the set of faces of a plane graph $G$. A face $f$ of length $k$ is called a $k$-\emph{face} and its length is denoted $|f|$. We write $(\geq k)$-\emph{face} and $(\leq k)$-\emph{face} to denote a face of length $\geq k$ and $\leq k$, respectively. A $4$-valent vertex is called an $(f_{_1},f_{_2},f_{_3},f_{_4})$-\emph{vertex}, if the faces incident with $v$ are of length $f_{_i}$, $1\leq i \leq 4$, and these are met in a clockwise order around $v$.   
 
\section{Linkages and wheels}\label{link} Throughout this section, $G$ is a $4$-connected plane graph, and  
\begin{equation}\label{uH}
\mbox{$u\in V(G) \sm V(X_{_G})$, $S_{_u} \subseteq H$ is the facial wheel of $u$, where $H$ is a $4$-hammock of $G$.}
\end{equation} 

By a $C_{_u}$-\emph{linkage} we mean a $(bnd H,C_{_u})$-linkage in $H$; such clearly does not meet $u$, by planarity. By $end \P$ we refer to the end vertices on $C_{_u}$ of members of a $C_{_u}$-\emph{linkage} $\P$. For such a $\P$, put $\alpha(\P) = |end \P \cap V(C_{_u}) \cap N(u)|$. Also, if $end \P =\{a_{_1},a_{_2},a_{_3},a_{_4}\}$, then we always assume these appear in this clockwise order along $C_{_u}$ and denote by $P_{_i}$ the member of $\P$ meeting $a_{_i}$. 

By planarity and since $N(u) \subseteq V(C_{_u})$, every $S_{_u} \cup \P$-bridge does not meet or attach to $u$. Let $P \in \P$ and let $P'$ be a member of $\P$ or a segment of $C_{_u}$. By $P$-\emph{ear} we mean an $S_{_u} \cup \P$-ear with both its ends in $P$. By $(P,P')$-\emph{ear} we mean an $S_{_u} \cup \P$-ear with one end in $P$ and the other in $P'$.

If for any $b \in (a_{_i} C_{_u} a_{_{i+1}})$ there exists a $C_{_u}$-linkage $\P'$ 
satisfying $end \P' = (end \P \sm \{a_{_i}\}) \cup \{b\}$ or $end \P' = (end \P \sm \{a_{_{i+1}}\}) \cup \{b\}$, then we call $\P$ \emph{slippery with respect to 
$[a_{_i} C_{_u} a_{_{i+1}}]$}, where $1 \leq i \leq 4$, and $a_{_5} = a_{_1}$.
We say that $\P$ is \emph{slippery} if it is slippery with respect to each segment 
$[a_{_i} C_{_u} a_{_{i+1}}]$ satisfying $a_{_i} C_{_u} a_{_{i+1}} \not= \emptyset$. 

\begin{statement}\label{slip}
A $C_{_u}$-linkage is slippery.
\end{statement}

\noindent
{\sl Proof.} Let $\P$ denote such a linkage, and let $w \in (a_{_i} C_{_u} a_{_{i+1}})$ such that $1 \leq i \leq 4$. Planarity and $C_{_u}$ being induced assert that there is an $S_{_u} \cup \P$-bridge $B$ with $w$ as an attachment. Such a bridge attaches to at least one of $P_{_i}-a_{_i}$ or $P_{_{i+1}} - a_{_{i+1}}$. This is clearly true if $B$ is trivial, as $C_{_u}$ is induced. If nontrivial,
then having all attachments of $B$ in $[a_{_i} C_{_u} a_{_{i+1}}]$ implies that 
the $3$-set consisting of $u$ and the two extremal attachments of $B$ on $[a_{_i} C_{_u} a_{_{i+1}}]$ is a $3$-disconnector of $G$, by planarity. \QED \\

It follows now from \sref{slip} that:

\begin{statement}\label{alpha>=1}
A $C_{_u}$-linkage satisfying $\alpha \geq 1$ exists. 
\end{statement} 

Our main tool for proving subsequent claims is the following. 
\begin{statement}\label{increase}
Suppose that: \\
(\sref{increase}.a) $H$ is a minimal fat $4$-hammock; and\\
(\sref{increase}.b) $\P$ is a $C_{_u}$-linkage with $end \P =\{a_{_1},a_{_2},a_{_3},a_{_4}\}$ satisfying:

\indent (\sref{increase}.b.1) $\alpha(\P) =k>0$, $k$ an integer; and

\indent (\sref{increase}.b.2) $a_{_1} ,a_{_3} \notin N(u)$; and 

\indent (\sref{increase}.b.3) $a_{_2} \in N(u)$; and

\indent (\sref{increase}.b.4) $|N(u) \cap [a_{_1}C_{_u} a_{_2}]|  \geq 2$. 

\noindent
Then, $K^-_{_4} \subseteq G$ or there exists a $C_{_u}$-linkage satisfying $\alpha \geq k+1$.
\end{statement}

\noindent
{\sl Proof.} Assume towards contradiction that 
\begin{equation}\label{assume}
\mbox{a $C_{_u}$-linkage with $\alpha \geq k+1$ does not exist.} 
\end{equation}

Let $P$ be the $a_{_1} a_{_3}$-segment of $C_{_u}$ not containing $a_{_4}$. By (\ref{assume}) and planarity, for any $\P$ satisfying (\sref{increase}.b), every member of $(N(u) \sm \{a_{_2}\}) \cap P$ is an attachment of an $S_{_u} \cup \P$-bridge attaching to $C_{_u}$ and $P_{_2}$ only. 
Such bridges exist by (\sref{increase}.b.4), (\ref{assume}), and since $C_{_u}$ is induced. Consequently, $a_{_2} \notin bnd H$. 

Choose a $\P$ satisfying (\sref{increase}.b) such that 
\begin{equation}\label{chooseP2}
\mbox{no $P_{_2}$-ears are embedded in the region of the plane interior to $[a_{_2} C_{_u} a_{_3}] \cup P_{_2} \cup P_{_3}$.}
\end{equation}

By (\sref{increase}.b.4), let $z \in N(u) \cap (a_{_1} C_{_u} a_{_2})$ such that $[a_{_1} P z]$ is minimal. Let $B$ be an $S_{_u} \cup \P$-bridge attached to $z$; such is embedded in the region of the plane interior to $[a_{_1} C_{_u} a_{_2}] \cup P_{_1} \cup P_{_2}$. 
By (\ref{assume}), 
\begin{equation}\label{notP1}
\mbox{$B$ has no attachment on $P_{_1}$.}
\end{equation}
Connectivity and existence of $z$ then imply that there are vertices $x \in [a_{_1} C_{_u} a_{_2}]$ (possibly $x=z$) and $y \in V(P_{_2})$ attachments of $B$ such that $[a_{_1}Px]$ and $[yP_{_2}v]$ are minimal, where $v \in V(P_{_2}) \cap bnd H$. 

By (\ref{assume}), 
\begin{equation}\label{notP2P3}
\mbox{there are no $(P_{_2}, P_{_3}-a_{_3})$-ears with an end in $[a_{_2} P_{_2} y)$.}
\end{equation}
Indeed, if such an ear exists, then $P_{_2}$ can be rerouted through $y$ and $B$ to meet $z$, and $P_{_3}$ can be rerouted through the ear and $[a_{_2}P_{_2}y)$ to meet $a_{_2}$; contradicting (\ref{assume}). 

Let $\ell \in [a_{_2} C_{_u} a_{_3}]$ be defined as follows. If there exist an $((a_{_2} P_{_2} y), [a_{_2} C_{_u} a_{_3}])$-ear, then $\ell$ is an end of such an ear such that $[\ell P a_{_3}]$ is minimal.
Otherwise, $\ell =a_{_2}$. 

By planarity, (\ref{chooseP2}), (\ref{notP1}), and (\ref{notP2P3}), $\{u,x,y,\ell\}$ form the boundary of a $4$-hammock of $H$; such is trivial or degenerate, by minimality of $H$. In either case, $x$ coincides with $z$ and $B$ consists of the single edge $xy$ (otherwise, there is a $k$-disconnector, $k \leq 3$, separating $B$ from the rest of $G$) implying that $\{x,u,a_{_2},y\}$ induce a $K^-_{_4}$. \QED \\

We infer the following from \sref{increase}.
\begin{statement}\label{consecutive}
Suppose $H$ is a minimal fat $4$-hammock. 
Then, $K^-_{_4}\subseteq G$ or there is a $C_{_u}$-linkage $\P$ satisfying:

(\sref{consecutive}.a) $\alpha(\P) \geq 2$; and

(\sref{consecutive}.b) if $\alpha \leq 2$ for every $C_{_u}$-linkage, then every $C_{_u}$-linkage with $\alpha =2$ meets $N(u)$ at\\ \indent \indent \indent \  consecutive members of $end \P$.
\end{statement}

\noindent
{\sl Proof.} A $C_{_u}$-linkage satisfying $\alpha \geq 1$ exists, by \sref{alpha>=1}. To show that such a linkage with $\alpha \geq 2$ exists, assume, towards contradiction, that every $C_{_u}$-linkage has $\alpha \leq 1$. Let $\P$ be a $C_{_u}$-linkage with $\alpha(\P) = 1$ and $end \P = \{a_{_1},a_{_2},a_{_3},a_{_4}\}$; choose such notation so that $a_{_2} \in N(u)$. As, by assumption, a linkage with $\alpha \geq 2$ does not exist, each vertex in $N(u) \sm \{a_{_2}\}$ is an attachment vertex of an $S_{_u} \cup \P$-bridge that has attachments on $C_{_u}$ and $P_{_2}$ only. 
Since $d(u) \geq 4$ and $C_{_u}$ is induced, such bridges exist and thus $a_{_2} \notin bnd H$. By rerouting $P_{_2}$ through such bridges we may choose such a $\P$ such that $N(u) \subseteq (a_{_1} C_{_u} a_{_2}]$. Thus, by \sref{increase} the claim follows. 

Suppose next, that $\alpha \leq 2$ for every $C_{_u}$-linkage, and suppose $\P$ is such a linkage with $\alpha(\P) =2$ so that $N(u)$ is met by nonconsecutive members of $end \P$, say, $a_{_2},a_{_4}$. As, by assumption, there is no linkage with $\alpha > 2$,  each vertex in $N(u) \sm \{a_{_2},a_{_4}\}$ is an attachment vertex of an $S_{_u} \cup \P$-bridge that has attachments on $C_{_u}$ and $P_{_2}$ only, or on $C_{_u}$ and $P_{_4}$ only (both options do not occur together). Since $d(u) \geq 4$ and $C_{_u}$ is induced, such bridges exist; hence $|bnd H \cap \{a_{_2}, a_{_4}\}| \leq 1$. By rerouting $P_{_2}$ and/or $P_{_4}$ through such bridges, we may choose $\P$ so that $|N(u) \cap (a_{_1} C_{_u} a_{_2}]| \geq 2$ or $|N(u) \cap [a_{_4} C_{_u} a_{_1})| \geq 2$. The claim then follows by \sref{increase}.\QED \\

We conclude this section with the following. 
\begin{statement}\label{short}
Let $H$ be minimal and fat and suppose $S_{_u}$ is short such that if it is imbalanced then it is proper. 
Then, a $C_{_u}$-linkage satisfying $\alpha \geq 3$ exists. 
\end{statement} 

\noindent
{\sl Proof.} Assume, to the contrary, that 
\begin{equation}\label{con}
\mbox{a $C_{_u}$-linkage satisfying $\alpha \geq 3$ does not exist.}
\end{equation} 
By \sref{consecutive}, a linkage with $\alpha =2$ exists; moreover, any $C_{_u}$-linkage  satisfying $\alpha=2$ meets $N(u)$ at consecutive ends. 
Suppose $\P$ is such a linkage where $end\P =\{a_{_1},a_{_2},a_{_3},a_{_4}\}$, and choose the notation so that the members of $end \P$ meeting $N(u)$ are $a_{_2}$ and $a_{_3}$. 

Since $C_{_u}$ is induced, 
\begin{equation}\label{meet}
(a_{_4} C_{_u} a_{_1}) \cap N(u) = \emptyset.
\end{equation}
Indeed, otherwise, a bridge attached to a member of $(a_{_4} C_{_u} a_{_1}) \cap N(u)$ has an attachment on  at least one of $P_{_1}$ or $P_{_4}$, by planarity and $4$-connectivity (see argument of \sref{slip}); contradicting (\ref{con}). 

Let $T,T'$ be as in (SH.1). $S_{_u}$ being short and \raf{meet} imply that 
either 
\begin{equation}\label{=1}
|V(T) \cap \{a_{_2},a_{_3}\}| = |V(T') \cap \{a_{_2},a_{_3}\}|=1,
\end{equation}
or
\begin{equation}\label{=1'}
|V(T'') \cap \{a_{_2},a_{_3}\}| = 2, T'' \in \{T,T'\}.
\end{equation}
In either case, \raf{meet} implies that 
$\{a_{_1},a_{_4}\} \subseteq int Q''$, $Q'' \in \{Q,Q'\}$, where $Q,Q'$ are as in (SH.2). Consequently, $S_{_u}$ is imbalanced and consequently proper, by assumption. That is, $bnd H \cap \{a_{_1},a_{_4}\} = \emptyset$. 

An $(S_{_u} - a_{_1},bnd H)$-linkage $\P'$ exists in $H-a_{_1}$; otherwise $a_{_1}$ and a $k$-disconnector, $k \leq 3$, separating $bnd H$ and $S_{_u} - a_{_1}$ in $H-a_{_1}$ form a proper $4$-hammock of $H$ that is neither trivial nor degenerate; contradicting the minimality of $H$. As, by assumption, $d(u) = 4$, $\P'$ does not meet $u$ and is a $C_{_u}$-linkage in $H$.

Since $S_{_u}$ is short, $\alpha(\P')=|end \P' \cap N(u)| \geq 2$. 
We may assume equality holds or the claim follows. 
If $N(u)$ is met by consecutive members of $\P'$, then these are not contained in a single triangle $T$ or $T'$, as this would contradict \raf{meet} (which applies to any $C_{_u}$-linkage with $\alpha =2$ meeting $N(u)$ at consecutive members). 
On the other hand, if $N(u)$ is met by nonconsecutive members of $\P'$ (so that \raf{=1} is satisfied by $\P'$), then a $C_{_u}$-linkage satisfying the premise of \sref{increase} exists (see argument of \sref{consecutive}) and the claim follows by \sref{increase}.\QED

\section{Short wheels in minimal fat hammocks} The purpose of this section is to prove \sref{main-dis}. Let $H$ be a minimal fat $4$-hammock of a $4$-connected plane graph $G$; such is $2$-connected, by \sref{2-con-ham}. Consequently, every member of $F(H)$ is a circuit of $H$, each edge of $H$ is contained in precisely $2$ faces (we use this in the proof of \raf{total} below), and each $v \in V(H)$ is incident with $d_H(v)$ distinct faces.   
A vertex $v \in V(H)$ is called \emph{good} if $d_H(v) \geq 5$ or $v \in bnd H$. 

\begin{statement}\label{main-dis}
Let $H$ be a minimal fat $4$-hammock of a $4$-connected plane graph $G$ satisfying:

(\sref{main-dis}.a) $K^-_{_4} \not\subseteq G$; and

(\sref{main-dis}.b) every $P_{_3} \iso P \subset H$ contains a good vertex; and

(\sref{main-dis}.c) every $K_{_3} \iso K \subset H$ contains $\geq 2$ good vertices.\\
\noindent
Then, $H$ contains a short facial wheel $S_{_u}$ for some $u \in V(H) \sm V(X_{_H})$ such that if $S_{_u}$ is imbalanced, then it is proper. 
\end{statement}

We shall use the well-known ``discharging method'' in order to prove \sref{main-dis}.
Such a method involves four main steps: (i) distributing initial charges to elements of the graph, (ii) calculating the total charge distributed using Euler's formula, (iii) redistributing  charges according to a set of discharging rules, and finally (iv) estimating the resultant charge of each element. In our case, we shall employ the following charging-discharging schemes.\\

\noindent
\scaps{Charging scheme.} For $x \in V(H) \cup F(H)$, define the \emph{charge} $ch(x)$ as follows:

(CH.1) $ch(v)=6-d_H(v)$, for any $v \in V(H)$.

(CH.2) $ch(f)=6-2|f|$, for any $f \in F(H) \sm \{X_{_H}\}$.

(CH.3) $ch(X_{_H})=-5\frac{2}{3}-2|X_{_H}|$. \\

Next, we show that

\begin{equation}\label{total}
\displaystyle\sum_{x \in V(H) \cup F(H)} ch(x) = \frac{1}{3}.
\end{equation}

\noindent
{\sl Proof.} 
\begin{eqnarray*}
\displaystyle\sum_{x \in V(H) \cup F(H)} ch(x) &=& -5\frac{2}{3} -2|X_H|+  \displaystyle\sum_{f \in F(H)\setminus X_H} (6-2|f|)  + \displaystyle\sum_{v \in V(H)} (6-d(v))\\ &=& -5\frac{2}{3} -2|X_H| + 6(|f(H)|-1) + \displaystyle\sum_{f \in F(H)\setminus X_H}(-2|f|) + \displaystyle\sum_{v \in V(H)} (6-d(v))\\
                &=& -5\frac{2}{3}+6(|f(H)|-1)-2(2|E|)+6|V(H)|-2|E(H)|\\
                &=& 6(F(H)-E(H)+V(H))-11\frac{2}{3}=\frac{1}{3}
\end{eqnarray*}
\QED \\

\noindent
\scaps{Discharging scheme.} In what follows, by \emph{send} we mean ``discharge'' or ``pass charge''. 
\begin{enumerate}
	\item[(DIS.1)] Let $v \in V(X_{_H})$, such that $2 \leq d_H(v) \leq 4$. 
		\begin{itemize}
			\item[(DIS.1.a)] 	If $d_H(v)=2$, then let $g$ the face incident with $v$ other 
			                  than $X_H$. 
			                  If $|g|=3$, then $v$ sends $4$ to $X_{_H}$. Otherwise $v$  
			                  sends $3\frac{2}{3}$ to $X_H$ and $\frac{1}{3}$ to $g$.
			                   
	    \item [(DIS.1.b)] If $d_H(v)=3$, then $v$ sends $2\frac{2}{3}$ to $X_{_H}$ and 
	                      $\frac{1}{3}$  to every incident $(\geq 4)$-face.   
	    \item  [(DIS.1.c)] If $d_H(v)=4$, then $v$ sends $1\frac{2}{3}$ to $X_{_H}$ and 
	                       $\frac{1}{3}$  to every incident $(\geq 4)$-face.
    \end{itemize}

	\item[(DIS.2)] If $v \in V(H)$ is at least $5$-valent, then, $v$ sends 
	               $\frac{1}{3}$ to every incident $(\geq 4)$-face. 
	\item[(DIS.3)] If $v\in V(H) \setminus V(X_H)$ is $4$-valent, then:
     \begin{enumerate}
     		\item[(DIS.3.a)] $v$ sends $\frac{2}{3}$ to every incident $4$-face.
	      \item[(DIS.3.b)] $v$ sends $1$ to every incident $5$-face, unless $v$ is a 
	                       $(3,4,3,5)$-vertex, and then $v$ sends $1\frac{1}{3}$ to its 
	                       single incident $5$-face.
	      \item [(DIS.3.c)] $v$ sends $1\frac{1}{3}$ to every $(\geq 6)$-face.
	   \end{enumerate}
\end{enumerate}

\noindent
\bfm{Proof of \sref{main-dis}.} Assume, to the contrary, that the claim is false and apply (CH.1-3) and (DIS.1-3) to members of $V(H) \cup F(H)$. Let $ch^*(x)$ denote the charge of a member of $V(H) \cup F(H)$ after applying (DIS.1-3). We obtain a contradiction to \raf{total} by showing that $ch^*(x) \leq 0$ for every $x \in V(H) \cup F(H)$. This is clearly implied by the following claims proved below.\\

\noindent
(\sref{main-dis}.A) \emph{$ch^*(v) \leq 0$, for each $v \in V(H)$.}\\
(\sref{main-dis}.B) \emph{$ch^*(f) \leq 0$, for each $f \in F(H) \sm \{X_{_H}\}$.}\\
(\sref{main-dis}.C) \emph{$ch^*(X_{_H}) \leq 0$.}\\

Observe that according to (DIS.1-3), faces do not send charge and vertices do not receive charge.\\

\noindent
{\sl Proof of} (\sref{main-dis}.A). It is sufficient to consider vertices $v$ satisfying 
$2 \leq d_H(v) \leq 4$. Indeed, if $d_H(v) \geq 6$, then $ch(v)=ch^*(v) \leq 0$ by (CH.1); and, if $d_H(v) = 5$, then $v$ is incident with at least three $(\geq 4)$-faces, as $K^-_{_4} \not\subseteq G$, implying that $ch^*(v) \leq 0$ by (DIS.2). 

By (DIS.1.a-c), $ch^*(v) \leq 0$ for every $v \in V(X_{_H})$ with $2\leq d_H(v) \leq4$. This is clear if $v$ is $2$-valent; and true in case $v$ is at least $3$-valent as such a vertex is incident with at least one $(\geq 4)$-face distinct of $X_{_H}$, since $K^-_{_4} \not\subseteq G$. 

It remains to consider $v \notin V(X_{_H})$ satisfying $2 \leq d_H(v) \leq 4$; such is clearly $4$-valent, as $\kappa(G) \geq 4$. We may assume $v$ is not incident with at least three $(\geq 4)$-faces, for otherwise $ch^*(v) \leq 0$ since by (DIS.3.a-c), $v$ sends at least $\frac{2}{3}$ to each $(\geq 4)$-face. Consequently, since $K^-_{_4} \not\subseteq G$, $v$ is incident with precisely two $3$-faces that are edge disjoint. Next, at least one of the remaining faces incident with $v$, say $f$, is a $4$-face for otherwise $ch^*(v) \leq 0$ by (DIS.3.b-c). The remaining face incident with $v$, say $g$, is a $(\geq 5)$-face for otherwise $H$ contains a short facial wheel; contradictory to our assumption. By (DIS.3.a), $v$ sends $2/3$ to $f$. Hence, $ch^*(v) \leq 0$ by (DIS.3.b) if $|g| = 5$, and by (DIS.3.c)  if $|g| \geq 6$.\inQED \\

\noindent
{\sl Proof of} (\sref{main-dis}.B). If $|f|=3$, then, $ch(f) = ch^*(f) = 0$ for any $3$-face $f$, by (CH.2). It remains to consider $(\geq 4)$-faces. If $f$ is such a face, then put $A_{_f} = \{v \in V(f) \sm bnd H: d_H(v) =4\}$ and note that (\sref{main-dis}.b) implies:
\begin{equation}
\label{eq:4}
\mbox{ $|A_f| \leq |f|-2$.}
\end{equation}
Clearly, 
\begin{equation}
\label{eq:main}
ch^*(f) = ch(f) + c(A_{_f}) + c(V(f)\sm A_{_f}),
\end{equation}   
where $c(X)$, $X \subseteq V(f)$, is the total charge sent to $f$ from members of $X$. 

We may assume that $f$ is a $5$-face. Indeed, if $|f|=4$, then $c(A_f)\leq \frac{2}{3}|A_f|$, by (DIS.1.c) and (DIS.3.a), $c(V(f) \setminus A_f)\leq \frac{1}{3}(|f|-|A_f|)$, by (DIS.1.a) and (DIS.2), and $|A_f| \leq 2$, by (\ref{eq:4}). Thus, $ch^*(f) \leq 0$, by \raf{eq:main}.
Next, if  $|f| \geq 6$, then $c(A_f)\leq 1\frac{1}{3}|A_f|$, by (DIS.1.c) and (DIS.3.c), $c(V(f) \setminus A_f)\leq \frac{1}{3}(|f|-|A_f|)$, by  (DIS.1.a) and (DIS.2), and $|A_f| \leq 4$, by (\ref{eq:4}). Hence, $ch^*(f) \leq 0$, by \raf{eq:main}. 

Assume then that $|f|=5$ so that $|A_f| \leq 3$, by \raf{eq:4}. We may assume that $f$ is incident with a $(3,4,3,5)$-vertex not in $V(X_H)$; otherwise, $c(A_f)\leq 1\times |A_f|$, by (DIS.1.c) and (DIS.3.b), $c(V(f)\setminus A_f) = \frac{1}{3}(|f|-|A_f|)$, by (DIS.1-2). By \raf{eq:main} (and as $|A_f| \leq 3$), $ch^*(f) \leq 0$. 

Let then $v \in V(f) \sm V(X_{_H})$ be a $(3,4,3,5)$-vertex. The members of $V(f)$ adjacent to $v$, say $v',v''$, are good by (\sref{main-dis}.c); and $|(V(f)\sm \{v,v',v''\}) \cap bndH|\geq 1$ or $S_{_u}$ is proper contradicting the assumption that such wheels do not exist in $H$. Let $v'''\in(V(f)\sm \{v,v',v''\}) \cap bndH$. $v$ sends $1\frac{1}{3}$ to $f$, By (DIS.3.c). Each of $\{v',v'',v'''\}$ sends $\frac{1}{3}$ to $f$, by (DIS.1-2) and since $f \neq X_{_H}$. The remaining vertex $V(f) \sm \{v,v',v'',v'''\}$ sends at most $1\frac{1}{3}$ to $f$, by (DIS.1-3) and since $f \neq X_H$. Consequently, $ch^*(f)=ch(f)+2\times 1\frac{1}{3} + 2\times \frac{1}{3}\leq 0$ (as $ch(f)=-4$).\inQED \\	

\noindent
{\sl Proof of} (\sref{main-dis}.C). For $i=2,\dots,5$, let
$A_i= \{v \in V(X_H):d_H(v)=i\}$;
$B=\{v \in V(X_H):d_H(v)\geq 5\}$;
$A_2'=\{v \in A_2:$ $v$ is incident with a 3-face$\}$; and put 
$A_2''=A_2 \setminus A_2'$.
Clearly, $A_i \subseteq bnd H$ for $i < 4$. Hence, since $H$ is a $4$-hammock of $G$ and $\kappa(G) \geq 4$,
\begin{equation}
\label{eq:6}
\mbox{$|A_2|+|A_3| \leq 4$.}
\end{equation}
By definition, 
\begin{equation}
\label{eq:7}
\mbox{$|A_2|+|A_3| +|A_4| \leq |X_H|$.}
\end{equation}

By (CH.3) and (DIS.1-2),
\begin{equation}
\label{eq:5}
\mbox{$ch^*(X_H)=-5\frac{2}{3}-2|X_H|+4|A_2'| +3\frac{2}{3}|A_2''|+2\frac{2}{3}|A_3|+1\frac{2}{3}|A_4|+\frac{1}{3}|A_5|$}
\end{equation}

By \raf{eq:5}, \raf{eq:6}, and \raf{eq:7}, it can be easily verified  that $ch^*(X_{_H}) \leq 0$ in the following cases:
(i) $|X_H| \geq 11$; (ii) $7 \leq |X_H| \leq 10$ and $|A_2|\not= 4$; and (iii) $4 \leq |X_H| \leq 6$ and $|A_2| \leq 2$. 

It remains to show that $ch^*(X_H) \leq 0$ in the cases: (I) $7 \leq |X_H| \leq 10$ and $|A_2|=4$ and (II) $4 \leq |X_H| \leq 6$, and $|A_2|\geq 3$. In the latter case, $V(X_H) \sm A_2$  is  a $k$-disconnector, $k \leq 3$, of $G$; this is so since $V(H) \sm V(X_H) \not= \emptyset$ by the fatness of $H$ and each vertex in $int H$ being at least $4$-valent. 

Suppose then that (I) occurs. Then, $|B|=0$ can be assumed; indeed, if $|B| \geq 1$, then $|A_2| +|A_3| +|A_4| \leq |X_H|-1$ implying that $ch^*(f)\leq 0$, by \raf{eq:5} and \raf{eq:6}. We may also assume that $|A_2'| \geq 1$; otherwise $|A_2''|=|A_2|=4$, and $ch^*(f)\leq 0$, by (\ref{eq:5}) and (\ref{eq:6}). Let then $x \in A_2' \subseteq bnd H$. $\{x\} \cup N_H(x)$ induce a $3$-face implying that at least one member of $N_H(x)$ is a good verex, by (\sref{main-dis}.c), and consequently in $bnd H$ as $A_5 \subseteq B = \emptyset$ (see above). 
As  $|X_H| \geq 7$ and thus $|V(H) \setminus \{x\}|\geq 6$, it follows that $(bnd H \sm \{x\}) \cup N_H(x)$ is either a $3$-disconnector of $G$ or a $4$-hammock of $H$ with its interior containing at least $2$ vertices; contradicting $\kappa(G) \geq 4$ and $H$ being minimal, respectively.\inQED \\
\QED

\section{Proof of \sref{main2}}\label{proof} Suppose $K^-_{_4} \not\subseteq G$ and let $v$ be an apex vertex of $G$ contained in some $5$-disconnector of $G$. Fix an embedding of $G$ and identify $G$ with its embedding. 

By~\cite[Lemma 2]{Mohar}(see Introduction), $\delta(G-v)=4$; implying that we may assume that $G-v$ has a minimal fat $4$-hammock $H$. To see this, let $u \in V(G-v)$ be $4$-valent. $N_{G-v}(u)$ is the boundary of two $4$-hammocks of $G-v$. If each of these two hammocks is degenerate, then $G$ is a $7$-vertex graph which contains a $TK_{_5}$. Thus, we may assume that at least one of these hammocks is fat; implying that minimal fat $4$-hammocks exist in $G-v$.

$H$ satisfies (\sref{main-dis}.b-c) or $K^-_{_4} \subseteq G$; hence, 
by \sref{main-dis}, there is a short facial wheel $S_{_u} \subseteq H$ with some $4$-valent vertex $u \notin V(X_{_H})$ as a hub; and such that $S_{_u}$ is proper if it is imbalanced.
Let $\P$ be a $C_{_u}$-linkage in $H$ satisfying $\alpha(\P) \geq 3$, by \sref{short}.
The set $\{v\} \cup bnd H$ forms the boundary of a $5$-hammock $H'$ of $G$ satisfying $S_{_u} \subseteq H'$; let $w \notin V(H')$ and let $F$ be a $(w, bnd H')$-$5$-fan in $G$, such clearly does not meet $int H'$. Observing that $uv \in E(G)$, as $u$ is $4$-valent in $G-v$, it follows that $TK_{_5} \subseteq S_{_u} \cup \P \cup F \cup \{uv\} \subseteq G$. \QED

\end{document}